\documentclass{article}
\usepackage[utf8]{inputenc}
\usepackage{mathtools} 
\usepackage{amsfonts} 
\usepackage{amsmath}
\usepackage{siunitx}
\usepackage{float}
\usepackage{listings}
\usepackage{nccmath}
\usepackage[english]{babel}
\usepackage{fancyhdr}

\usepackage{listings,xcolor}

\title{On a variant of the happy numbers and their generalizations}
\author{Luca Onnis}
\date{February 2022}
\begin{document}

\maketitle
\begin{abstract}
    This paper investigates a variant of the famous "happy numbers" sequence, given by A351327 \cite{oeis2} on the oeis. First of all we'll define this integer sequence, and then we'll show some important results about it; in particular we conjectured that if $k$ is a term of the sequence, then it converges to 1 in a number of steps less or equal to 3. Furthermore it will be possible to find some codes written in Wolfram language in order to compute large terms of the sequence and to support our hypothesis. At the end we'll explore some new conjectures and generalizations about this kind of integer sequence.
\end{abstract}
\section{Introduction}
\subsection{Happy numbers}
The sequence is given by A007770 \cite{oeis1} on the oeis. \\
To determine whether a given number $k$ is a term of this sequence, start with $k$, take the square of each digit and sum them together, apply the same process to the result, and continue until 1 is reached or a loop is entered. If 1 is reached, $k$ is a term of this sequence. If not, $k$ is a unhappy number. For example, consider the number 7: \\
7 is happy because after 5 iterations we reach the number 1. In fact:
\[
7\rightarrow49\rightarrow97\rightarrow130\rightarrow10\rightarrow1
\]
\subsection{Variant of happy numbers}
To determine whether a given number k is a term of this sequence, start with k, take the square of the product of its nonzero digits, apply the same process to the result, and continue until 1 is reached or a loop is entered. If 1 is reached, k is a term of this sequence. \\
For example, consider the number 375. 375 is a term of the sequence, in fact:
\[
375\rightarrow11025\rightarrow100\rightarrow1
\]
Now consider the number 4. 4 is not a term, in fact:
\[
4\rightarrow16\rightarrow36\rightarrow324\rightarrow576\rightarrow44100\rightarrow256\rightarrow3600\rightarrow324
\]
and we reached a loop of length 5 (starting with 324).
\subsubsection{First observations regarding this variant}
First of all we can notice that, since 1 is the neutral element for moltiplication in $\mathbb{N}$, if $k$ is a term, the numbers obtained by inserting ones anywhere in $k$ are terms. \\
Furthermore, since we are multiplying nonzero digits, we also have that if $k$ is a term, the numbers obtained by inserting zeros anywhere in $k$ are terms. \\
We also know that moltiplication is commutative so if $k$ is a term, each distinct permutation of the digits of $k$ gives another term. \\
This sequence counts an infinite amount of terms, in fact every power of 10 is a term. In fact:
\[
10^{k}\rightarrow 1 \mbox{ for every positive integer $k$ }
\]
This is true because $10^{k}$ has a 1 followed by $k$ zeros, so multiplying the nonzero digits together, we'll obviously get 1 at first iteration. \\
In order to compute terms less or equal to $h$ you can use this code on Wolfram Mathematica:
\begin{lstlisting}[language=Mathematica,caption={To compute terms less or equal to $h$}]
Select[Range[h],
FixedPoint[
    Product[ReplaceAll[0 -> 1][IntegerDigits[#]][[i]]^2, {i, 1,
       Length[ReplaceAll[0 -> 1][IntegerDigits[#]]]}] &, #, 10] == 1 &]
  \end{lstlisting}
\section{Lemma 1}
Let $k$ be a term of the sequence, then the second last iteration is a number of the form $10^{2h}$ where $h\in\mathbb{N}$
\subsection{Proof}
We want to show that:
\[
k\rightarrow k_1\rightarrow k_2\rightarrow\dots\rightarrow 10^{2h}\rightarrow1
\]
First of all note that if 1 is reached, then the second last iteration is a number composed by only ones and zeros. Otherwise,we can assume that there is a digit strictly bigger than 1 , call it $a_0$; then the square of the product of the nonzero digits is at least $a_0^{2}$ which is strictly bigger than 1. \\
Consider a natural number composed by only ones and zeros:
\[
n_0=\sum_{b=0}^{m}a_b\cdot10^{b} \mbox{ where $0\leq a_b\leq 1$ } \mbox{ $\forall b\in\{0,1,\dots,m\}$}
\]
then consider a subset of $\{0,1,\dots,m\}$:
\[
\{u_1,u_2,\dots,u_p\}\subseteq\{0,1,\dots,m\} \mbox{ such that $a_{u_i}=1$} \mbox{ $\forall i\in \{1,\dots,p\}$}
\]
So we'll have that
\[
n_0=\sum_{i=1}^{p}a_{u_i}\cdot10^{u_i}=\sum_{i=1}^{p}10^{u_i}
\]
and:
\[
\Bigl[\prod_{i=1}^{p}a_{u_i}\Bigr]^{2}=1
\]
We want to show that the only 1 in the decimal representation of the number of the second last iteration is the first. \\
Now consider the third last iteration:
\[
\dots\rightarrow m\rightarrow  N\rightarrow1 \mbox{ where $n_0$ has more than a 1 in its decimal representation}
\]
Consider $P(n)$ as the function which compute the product of the nonzero digits of $n$. If such $N$ exists, then:
\[
[P(m)]^{2}=n_0
\]
We'll prove by contradiction that such number can not exist.\\
$N$ is a perfect square whose prime factor are all less than or equal to 7. That's because $P(m)$ is the product of the digits of $m$, and the digits in base 10 are 0,1,2,3,4,5,6,7,8,9 where the only prime numbers are 2,3,5,7. \\
We can assume that $N$ last digit is 1. If not, it ends with $2q$ zeros, so $N=2^{2q}\cdot5^{2q}\cdot n_1$ where $N_1$ is a number composed by ones and zeros whose last digit is 1. \\
Let $K=\{n\in\mathbb{N}$ such that $n$ is a perfect square whose decimal representation is composed only by zeros and ones$\}\setminus\{10^{f}:f\in\mathbb{N}\}$. \\
For the Well-ordering principle, it exists $\min(K)=N$. \\
If:
\[
d^{2}=N\equiv 0 \mod (10)
\]
then $d=10d_1$ , so $100\mid N$ and $(\frac{d}{10})^{2}$ is a perfect square strictly lower than $N$, with the required characteristics. This is in contradiction with our choice of $N$ (which is the minimum of the set $K$).
Otherwise $N$ ends with "01" or "11". \\
But if:
\[
d^{2}=N\equiv 11 \mod(100) \Rightarrow d^{2}\equiv 11 \equiv 3 \mod 4
\]
and there are no integer solution for this equation. \\
If:
\[
d^{2}=N\equiv 1 \mod(100)
\]
Then:
\[
d=100d_1+1 \mbox{ or } d=100d_1+49 \mbox{ or } d=100d_1+51 \mbox{ or } d=100d_1+99
\]
\[
\Rightarrow 2\not| d \wedge 5\not| d
\]
But $d$ is a product of digits, so:
\[
3\mid d \mbox{ $\vee$ } 7\mid d
\]
And no other primes can divide $d$. Since $9\mid N$ and $N$ is composed only by zeros and ones, we can conclude that the number of ones in $N$ is a multiple of 9. \\
So $d^{2}=3^{2a}\cdot 7^{2b}=9^{a}\cdot 49^{b}=N$ for some $a,b\in\mathbb{N}$. \\
Note that $a\equiv b \mod (2)$ because $N\equiv 1 \mod (10)$. So if $a$ is even, $b$ is even too and if $a$ is odd, $b$ is odd too. Furthermore we know that $N\equiv 1 \mod(100)$; Note that:
\[
9\cdot49^{2r+1}\equiv 41 \mod(100) \mbox{ $\forall r\in\mathbb{N}$}
\]
\[
9^{2}\cdot49^{2r}\equiv 81 \mod(100) \mbox{ $\forall r\in\mathbb{N}$}
\]
\[
9^{3}\cdot49^{2r+1}\equiv 21 \mod(100) \mbox{ $\forall r\in\mathbb{N}$}
\]
\[
9^{4}\cdot49^{2r}\equiv 61 \mod(100) \mbox{ $\forall r\in\mathbb{N}$}
\]
\[
9^{5}\cdot49^{2r+1}\equiv 1 \mod(100) \mbox{ $\forall r\in\mathbb{N}$}
\]
\[
\vdots
\]
\[
9^{10}\cdot49^{2r}\equiv 1 \mod(100) \mbox{ $\forall r\in\mathbb{N}$}
\]
\[
9^{15}\cdot49^{2r+1}\equiv 1 \mod(100) \mbox{ $\forall r\in\mathbb{N}$}
\]
In conclusion:
\[
9^{10t+5}\cdot49^{2r+1}\equiv 1 \mod(100) \mbox{ $\forall r,t\in\mathbb{N}$}
\]
And
\[
9^{10t}\cdot49^{2r}\equiv 1 \mod(100) \mbox{ $\forall r,t\in\mathbb{N}$}
\]
Note that:
\[
n^{2}\equiv 1,4,0 \mod(8)
\]
so $N=d^{2}$ can not end with "101" because $101\equiv 5 \mod 8$, so $N\equiv 1 \mod(1000)$
It's possible to continue this process in order to exclude other numbers, although it remains an oper problem to understsand if such number can exists. Using an heuristic argument we conclude that the probability to find a perfect square number of at most $n$ digits composed only by ones and zeros is equal to:
\[
\frac{1}{(\sqrt{5}\cdot10)^{n}} \mbox{ \cite{math}}
\]
For example,using this code: 
\begin{lstlisting}[language=Mathematica,caption={}]
Table[Mod[9^10*49^(838938 + 2500000 n), 10^9], {n, 1, 15}]
\end{lstlisting}
we know that:
\[
9^{10}\cdot 49^{838938+2500000t}\equiv 111110001 \mod 10^9 \mbox{ $\forall t \in\mathbb{N}$}
\]
and for $t=1$ we have that:
\begin{lstlisting}[language=Mathematica,caption={}]
Length[IntegerDigits[9^10*49^(838938 + 2500000)]]
=5643470
\end{lstlisting}
so:
\[
9^{10}\cdot 49^{3338938} \mbox{ has 5643470 digits}
\]
\section{Conjecture 1}
There are no perfect square that are term of the sequence with length of convergence equal to 3. \\
Assuming this conjecture, it follows that if $n_0$ is a perfect square and a term of this sequence, then its process to converge to 1 is given by:
\[
n_0\rightarrow 1
\]
if $n_0=10^{2h}$ for some $h\in\mathbb{N}$. Or:
\[
n_0\rightarrow 10^{2h}\rightarrow 1
\]
if $n_0$ is a perfect square composed by zeros, ones, $h$ fives, $x$ twos, $y$ fours, $z$ eights such that $x+2y+3z=h$. This conditions are necessary because the product of the nonzero digits of $n_0$ , call it $P(n_0)$, is equal to:
\[
P(n_0)=10^{2h}=2^{h}\cdot 5^{h}=2^{x}\cdot 4^{y}\cdot 8^{z}\cdot 5^{h}= 2^{x+2y+3z}\cdot 5^{h}
\]
Using this code:
\begin{lstlisting}[language=Mathematica,caption={}]
Select[Range[1, 10000], 
 Complement[IntegerDigits[#^2], {1, 0, 2, 5, 4, 8}] == {} && 
   DigitCount[#^2, 10, 5] == 
    DigitCount[#^2, 10, 2] + 2 DigitCount[#^2, 10, 4] + 
     3 DigitCount[#^2, 10, 8] &]
  \end{lstlisting}
it's possible to compute numbers less or equal to 10000 such that their square follows the rules above. \\
For example, one of them is 7152 ($7152^{2}=51151104$ and $[P(51151104)]^{2}=10000$). So its process to converge to 1 is:
\[
7152^{2}\rightarrow 10^{4}\rightarrow 1
\]
But there aren't numbers whose product of their nonzero digits is equal to 7152. That's because $7152=2^{4}\cdot 3 \cdot 149$. In the last section there will be a deeper exploration of this fact.
\section{Theorem 1}
Let $k$ be a term of the sequence, then it converges to 1 in a number of steps less than or equal to 3.
\subsection{Proof of Theorem 1}
First of all note that theorem 1 is true if and only if the fourth last iteration is not a perfect square. \\
Suppose that Theorem 1 is true, so consider the following process of length 3:
\[
n_0\rightarrow n_1 \rightarrow 10^{2h} \rightarrow 1
\]
Assuming conjecture 1, if $n_0$ is a perfect square and is a term of the sequence then it converges to 1 in a number of steps less or equal to 2 while our $n_0$ converges to 1 in 3 steps by hypothesis. \\
Now suppose that the fourth last iteration is not a perfect square and that $n_0$ converges to 1 in a number of steps $k>3$. Then the process is of the form:
\[
n_0\rightarrow n_1 \rightarrow\dots\rightarrow n_{k-3}\rightarrow n_{k-2}\rightarrow n_{k-1} \rightarrow 1
\]
But the fourth last iteration should be the square of the nonzero digits of $n_{k-4}$, so this is a contradiction caused by supposing that $n_0$ converges to 1 in a number of steps strictly bigger than 3. \\
Again, assuming conjecture 1 is obvious that the fourth last iteration is not a perfect square. That's because all perfect squares that are terms of the sequence, converge to 1 in a number of steps less than or equal to 2. But the number in the fourth last iteration converge to 1 in $3>2$ steps.
\section{Other conjectures and generalizations}
In conjecture 1 I said that if you consider the process:
\[
n_0\rightarrow n_1 \rightarrow 10^{2h}\rightarrow 1
\]
Then $n_1$ must be composed by zeros, ones, $h$ fives, $x$ twos, $y$ fours, $z$ eights such that $x+2y+3z=h$. But testing out the first terms below $10^{8}$ that converge to 1 in 3 steps, and their process of convergence, I note that their second iteration is always a number composed by only zeros, ones, $h$ fives and $h$ twos. So there are no fours and no eights in their decimal representation. \\
In order to compute the entire process it's possible to use this code:
\begin{lstlisting}[language=Mathematica,caption={}]
Table[FixedPointList[
   Product[ReplaceAll[0 -> 1][IntegerDigits[#]][[i]]^2, {i, 1, 
      Length[ReplaceAll[0 -> 1][IntegerDigits[#]]]}] &, 
   Part[Select[Range[1000], 
     FixedPoint[
        Product[ReplaceAll[0 -> 1][IntegerDigits[#]][[i]]^2, {i, 1, 
           Length[ReplaceAll[0 -> 1][IntegerDigits[#]]]}] &, #, 10] ==
        1 &], j], 10], {j, 1, 
   Length[Select[Range[1000], 
     FixedPoint[
        Product[ReplaceAll[0 -> 1][IntegerDigits[#]][[i]]^2, {i, 1, 
           Length[ReplaceAll[0 -> 1][IntegerDigits[#]]]}] &, #, 10] ==
        1 &]]}] // TableForm
          \end{lstlisting}
          which gives for the first terms below 1000:
\begin{lstlisting}[language=Mathematica,caption={}]
1	1			
5	25	100	1	1
10	1	1		
11	1	1		
15	25	100	1	1
25	100	1	1	
50	25	100	1	1
51	25	100	1	1
52	100	1	1	
100	1	1		
101	1	1		
105	25	100	1	1
110	1	1		
111	1	1		
115	25	100	1	1
125	100	1	1	
150	25	100	1	1
151	25	100	1	1
152	100	1	1	
205	100	1	1	
215	100	1	1	
250	100	1	1	
251	100	1	1	
255	2500	100	1	1
357	11025	100	1	1
375	11025	100	1	1
455	10000	1	1	
500	25	100	1	1
501	25	100	1	1
502	100	1	1	
510	25	100	1	1
511	25	100	1	1
512	100	1	1	
520	100	1	1	
521	100	1	1	
525	2500	100	1	1
537	11025	100	1	1
545	10000	1	1	
552	2500	100	1	1
554	10000	1	1	
573	11025	100	1	1
735	11025	100	1	1
753	11025	100	1	1
1000	1	1		
        \end{lstlisting}
Using this conjecture we think that the only numbers, whose square is composed by zeros, ones, $h$ fives and $h$ twos, and whose prime factor are in the set $\{2,3,5,7\}$ are all of the form:
\[
10^{m} \mbox{ or } 5\cdot 10^{m} \mbox{ or } 105\cdot 10^{m} \mbox{ for some $m\in\mathbb{N}$}
\]
as suggested by using this code:
\begin{lstlisting}[language=Mathematica,caption={}]
Intersection[
 Select[Range[1, 1000000], 
  Complement[IntegerDigits[#^2], {1, 0, 2, 5}] == {} && 
    DigitCount[#^2, 10, 5] == DigitCount[#^2, 10, 2] &], 
 Select[Range[1, 1000000], 
  Complement[Part[FactorInteger[#], All, 1], {2, 3, 5, 7}] == {} &]]
 \end{lstlisting}
It's also possible to notice that if $k$ is a term of the sequence, then there's no 9 in its decimal representation (Verified for terms below $10^{10}$). \\
\subsection{Higher exponents}
Instead of squaring the digits, you can actually do the same process raising each nonzero digits with higher exponents than 2. A very interesting case is when you take the cube of the product of the nonzero digits during the process. \\
I define $S_k$ as the set of the terms of the sequence which follow the process raising the product of their nonzero digits to the power of $k$. I define $S_k(n)$ as the finite set that contains all terms less or equal to $n$. Then, it's possible to note that the cardinality of $S_2(n)$ is lower than the cardinality of $S_3(n)$. So:
\[
|S_2(n)|\leq|S_3(n)| \mbox{ $\forall n\in\mathbb{N}$}
\]
In fact:
\begin{lstlisting}[language=Mathematica,caption={To compute terms below 1000 of $S_2$}]
Select[Range[1000], 
 FixedPoint[
    Product[ReplaceAll[0 -> 1][IntegerDigits[#]][[i]]^2, {i, 1, 
       Length[ReplaceAll[0 -> 1][IntegerDigits[#]]]}] &, #, 10] == 
   1 &]
          \end{lstlisting}
          which gives:
 \begin{lstlisting}[language=Mathematica,caption={Terms below 1000 of $S_2$}]
1, 5, 10, 11, 15, 25, 50, 51, 52, 100, 101, 105, 110, 111, 115, 125, 
150, 151, 152, 205, 215, 250, 251, 255, 357, 375, 455, 500, 501, 502, 
510, 511, 512, 520, 521, 525, 537, 545, 552, 554, 573, 735, 753, 1000
          \end{lstlisting}     
While:
\begin{lstlisting}[language=Mathematica,caption={To compute terms below 1000 of $S_3$}]
Select[Range[1000], 
 FixedPoint[
    Product[ReplaceAll[0 -> 1][IntegerDigits[#]][[i]]^3, {i, 1, 
       Length[ReplaceAll[0 -> 1][IntegerDigits[#]]]}] &, #, 10] == 
   1 &]
          \end{lstlisting}
          which gives:
       \begin{lstlisting}[language=Mathematica,caption={Terms below 1000 of $S_3$}]
1, 2, 3, 5, 8, 10, 11, 12, 13, 15, 18, 20, 21, 24, 25, 27, 30, 31, 
42, 45, 50, 51, 52, 54, 55, 56, 57, 65, 72, 75, 80, 81, 100, 101, 
102, 103, 105, 108, 110, 111, 112, 113, 115, 118, 120, 121, 124, 125, 
127, 130, 131, 142, 145, 150, 151, 152, 154, 155, 156, 157, 165, 172, 
175, 180, 181, 200, 201, 204, 205, 207, 210, 211, 214, 215, 217, 222, 
225, 235, 240, 241, 250, 251, 252, 253, 255, 258, 270, 271, 285, 300, 
301, 310, 311, 325, 352, 355, 377, 402, 405, 412, 415, 420, 421, 445, 
450, 451, 454, 455, 457, 475, 478, 487, 500, 501, 502, 504, 505, 506, 
507, 510, 511, 512, 514, 515, 516, 517, 520, 521, 522, 523, 525, 528, 
532, 535, 540, 541, 544, 545, 547, 550, 551, 552, 553, 554, 558, 560, 
561, 570, 571, 574, 582, 585, 605, 615, 650, 651, 679, 697, 702, 705, 
712, 715, 720, 721, 737, 745, 748, 750, 751, 754, 769, 773, 784, 796, 
800, 801, 810, 811, 825, 847, 852, 855, 874, 967, 976, 1000
          \end{lstlisting}   
It is conjectured that if $k$ is a term of the sequence $S_3$, then it converges to 1 in a number of steps less than or equal to 10. Maybe this is the reason why there are more terms in $S_3(n)$ than in $S_2(n)$. Although there are some numbers that during the process seems to diverge to infinity. For example the number 4:
\[
4\rightarrow 64\rightarrow13824\rightarrow7077888\rightarrow5416169448144896\rightarrow
\]
\[
\rightarrow188436971148778297205194752000\rightarrow
\]
\[
\rightarrow1545896640285238037724131582088286996267008000000\rightarrow\dots
\]
Instead the number 217 converge to 1 in 8 steps:
\[
217\rightarrow2744\rightarrow11239424\rightarrow5159780352\rightarrow54010152000000000\rightarrow8000000\rightarrow512\rightarrow1000\rightarrow1
\]
For completeness I'll put the terms below 1000 of $S_k$ for every $2\leq k\leq 5$.
 \begin{lstlisting}[language=Mathematica,caption={Terms below 1000 of $S_4$}]
1, 10, 11, 25, 52, 100, 101, 110, 111, 125, 152, 205, 215, 250, 251, 
455, 502, 512, 520, 521, 545, 554, 1000
          \end{lstlisting}
\begin{lstlisting}[language=Mathematica,caption={Terms below 1000 of $S_5$}]
1, 10, 11, 25, 52, 100, 101, 110, 111, 125, 152, 205, 215, 250, 251, 
455, 502, 512, 520, 521, 545, 554, 1000
          \end{lstlisting}  
We think that using a sufficiently big $k_0$ we could conclude that all terms of $S_{k}$ for $k>k_0$ converge to 1 in a number of steps less or equal to 2.
  \bibliographystyle{plain}
\bibliography{Bibliografia.bib}
\end{document}